\documentclass[10pt, a4paper]{amsart}
\usepackage{latexsym, amsfonts, euscript}

\def\Z{{\mathbb Z}}

\def\PP{{\mathbb P}}

\def\Fq{{\mathbb F}_q}

\def\supp{{\rm supp}}

\def\a{{\alpha}}
\def\b{{\beta}}

\def\ilm{{I_{\ell}[m]}}

\def\L{{\Lambda}}

\newcommand\bprf{\begin{proof}}
\newcommand\eprf{\end{proof}}
\newtheorem{theorem}{Theorem}
\newtheorem{proposition}[theorem]{Proposition}
\newtheorem{lemma}[theorem]{Lemma}
\newtheorem{corollary}[theorem]{Corollary}
\newtheorem{conjecture}[theorem]{Conjecture}
\newtheorem*{theorem*}{Theorem}
\newtheorem*{conjecture*}{Conjecture}
\theoremstyle{remark}

\newtheorem{example}[theorem]{Example}

\begin{document}

\title{Subclose Families,  Threshold Graphs, and the 
    Weight Hierarchy of Grassmann and Schubert Codes} 

\author{Sudhir R. Ghorpade}
\address{Department of Mathematics, 
Indian Institute of Technology Bombay,\newline \indent
Powai, Mumbai 400076, India.}
\email{srg@math.iitb.ac.in}

\author{Arunkumar R. Patil}
\address{Department of Mathematics,
Indian Institute of Technology Bombay,\newline \indent
Powai, Mumbai 400076, India \newline \indent
and \newline \indent
Shri Guru Gobind Shinghji Institute  of Engineering \& Technology,\newline \indent
Vishnupuri, Nanded 431 606, India}
\email{arun.iitb@gmail.com}

\author{Harish K. Pillai}
\address{Department of Electrical Engineering,
Indian Institute of Technology Bombay,\newline \indent
Powai, Mumbai 400076, India.}
\email{hp@ee.iitb.ac.in}

\date{\today} 

\keywords{Linear code, higher weight,  Grassmann variety, Grassmann code, Schubert variety, Schubert code, threshold graph, optimal graph.}

\subjclass[2000]{94B05, 05C07, 05C35, 14M15}
 
\dedicatory{Dedicated to Gilles Lachaud on his sixtieth birthday}

\begin{abstract}
We discuss the problem of determining the complete weight hierarchy of 
linear error correcting codes associated to Grassmann varieties and, more generally, to Schubert varieties in Grassmannians.  In geometric terms, this corresponds to determining  the maximum number of $\Fq$-rational points on sections of Schubert varieties (with nondegenerate Pl\"ucker embedding) 
by linear subvarieties of a fixed (co)dimension. The problem is partially solved in the case of Grassmann codes, and one of the solutions uses 
the combinatorial notion of a closed family. We propose a generalization of this to what is called a subclose family. A number of properties of subclose families are proved, and its connection with the notion of threshold graphs and graphs with maximum sum of squares of vertex degrees is outlined. 
\end{abstract}

\maketitle


\section{Introduction}
\label{sec:intro}

It has been almost a decade since the first named author and Gilles Lachaud wrote \cite{GL} 
where alternative proofs of Nogin's results on higher weights of Grassmann codes  \cite{N} were given and Schubert codes were introduced. Originally, much of \cite{GL} was conceived as a side remark in \cite{GL2}. But in retrospect, it appears to have been a good idea to write \cite{GL} as an independent article and use the opportunity to propose therein a conjecture concerning the minimum distance of Schubert codes. This conjecture has been of some interest, and 
after being proved, in the affirmative, in a number of special cases (cf. \cite{chen,V,GT,GV}), the 
general case appears to have been settled very recently by Xiang \cite{xiang}. 
The time sees ripe, therefore, to up the ante and think about more general questions. 
It is with this in view, that we discuss in this paper the problem of determining the complete weight hierarchy of Schubert codes and, in particular, the Grassmann codes. In fact, the case of Grassmann codes and the determination of higher weights in the cases 
not covered by the result in \cite{N} and \cite{GL} has already been considered in some recent work (cf. \cite{HJR,HJR2,GPP}). 
What is proposed here is basically a plausible approach to tackle the general case. 
%

This paper is organized as follows. In Section \ref{sec1} below we recall the combinatorial notion of a close family, 
and introduce a more general notion of a subclose family. A number of elementary properties of subclose families are proved here, including a nice duality that prevails among these. 
Basic notions concerning linear error correcting codes, such as 
the minimum distance and more generally, the higher weights 
are reviewed in Section \ref{sec2}. Further, we state here a general conjecture that relates the higher weights of
Grassmann codes and Schubert codes with subclose families. 
Finally, in Section \ref{sec3}, we recall threshold graphs and optimal graphs and then show that, in a special 
case, subclose families are closely related to these well-studied notions in graph theory. As an application we obtain  explicit bounds on the sum of squares of degrees of a simple 
graph in terms of the number of vertices and edges, which  seem to ameliorate and complement some of the known results 
on this topic that has been of some interest in graph theory  (cf. \cite{caen,das,peled}).  

\section{Close Families and Subclose Families}
\label{sec1}

Fix integers $\ell,m$ such that $1\le \ell \le m$. Set 
$$
k:= {{m}\choose {\ell}} \quad \mbox{ and } \quad \mu:=\max\{\ell,m-\ell\} +1. 
$$
Let $[m]$ denote the set $\{1, \dots , m\}$ of first $m$ positive integers. Given any nonnegative integer $j$, let $I_j[m]$ denote the set of all subsets of $[m]$ of cardinality $j$. 

Let $\Lambda \subseteq \ilm$. Following \cite{GL2}, we call $\L$ a \emph{close family} if $|A\cap B| = \ell -1$ for all $A,B\in \Lambda$. Suppose $|\L|=r$. Then $\L$ is said to be of \emph{Type I} if if there exists $S\in I_{\ell-1}[m]$ and
$T\subseteq [m]\setminus S$ with $|T|=r$ such that 
$$
\L = \{ S\cup \{ t\} : t \in T \},
$$
whereas $\L$ is said to be of \emph{Type II} if there exists $S\in I_{\ell-r+1}[m]$ and
$T\subseteq [m]\setminus S$ with $|T|=r$ such that 
$$
\L = \{ S\cup T \setminus \{ t\} : t \in T \}.
$$
Basic results about close families are as follows. 

\begin{proposition}[Structure Theorem for Close Families]
\label{StrThm}
Let $\L \subseteq \ilm$. 
Then $\L$ is close if and only if $\L$ is either of Type I or of  Type II.
\end{proposition}

This is proved in \cite[Thm. 4.2]{GL2}. An immediate consequence is the following. 
\begin{corollary}
\label{cor:uptomu}
Let $r$ be a nonnegative integer. 
A close family in $\ilm$ of cardinality $r$ exists if and only if   $r \le \mu$. 
In greater details, a close family of Type I in $\ilm$ exists if and only if $r\le m-\ell +1$, whereas 
 a close family of Type II in $\ilm$ exists if and only if $r\le \ell +1$.
\end{corollary}

We use this opportunity to state the following elementary result which complements Proposition \ref{StrThm}. This is not stated explicitly in \cite{GL,GL2}, but a related result is proved in \cite{GPP} where we obtain an algebraic counterpart of Proposition \ref{StrThm} in the setting of exterior algebras and the Hodge star operator.

\begin{proposition}[Duality]
\label{duality}
Given $\L\subseteq \ilm$, let ${\L}^{\ast}:=\{[m]\setminus A : A\in \L\} \subseteq I_{m-\ell}[m]$. Then 
$\L$ is close in $\ilm$ of type I if and only if $\L^\ast$ is close in $ I_{m-\ell}[m]$ of type II. 
\end{proposition}

\bprf
Given $S\in I_{\ell-1}[m]$ and
$T\subseteq [m]\setminus S$ with $|T|=r$, observe that 
$$
[m]\setminus \left(S\cup\{t\}\right) = \left([m]\setminus (S\cup T)\right)\cup T \setminus \{t\} 
$$
for every  $t\in T$.
\eprf

As explained in \cite{GL}, Corollary \ref{cor:uptomu} essentially accounts for the 
barrier on $r$ 
for which the higher weights $d_r$ of Grassmann codes $C(\ell,m)$ are hitherto known. (See, e.g., \cite{N,GL}.)
Recently some attempts have been made to break this barrier (cf. \cite{HJR,HJR2,GPP}) but the complete weight hierarchy
$\{d_r : 1\le r\le k\}$ 
is still not known. We will comment more on this in Section \ref{sec2}. 
For the time being, we introduce a combinatorial generalization of close families which may play some role in 
the determination of higher weights. 

Given a subset $\Lambda=\left\{A_1,\cdots,A_r\right\}$ of $\ilm$, we define
$$
K_\Lambda=\displaystyle\sum_{i<j}|A_i\cap A_j|.
$$
Further, given any nonnegative integer $r\le k$, 
we define 
$$
K_r(\ell,m):=\max \{K_\Lambda : \Lambda \subseteq \ilm \mbox{ and } |\Lambda|=r\}.
$$
We call $\L$ a \emph{subclose family} if $K_{\L} = K_r(\ell,m)$ where $r=|\L|$. 
It is clear that for each nonnegative integer $r\le k$, 
there exists a subclose family 
of cardinality $r$.

\begin{proposition}
\label{subcloseisclose}
Given $\L\subseteq \ilm$, we have 
$$
K_{\L} \le (\ell -1) {{|\L|}\choose {2}}.
$$
Moreover, equality holds if and only if $\L$ is a close family.
Consequently, 
for any nonnegative integer $r\le k$, 
we have 
$$
K_r(\ell,m)\le (\ell -1) {{r}\choose {2}}. 
$$
Moreover, equality holds if and only if $r\le \mu$.
\end{proposition}

\bprf
The last assertion follows from Corollary \ref{cor:uptomu}. The remaining assertions are obvious. 
\eprf

It is an interesting question to determine $K_r(\ell,m)$ for any $r$. The first few values are given by the above result. We shall now determine some more. To this end, let us first observe the following analogue of 
Proposition \ref{duality}. 

\begin{proposition}[First Duality Theorem]
\label{firstduality}
Given $\L\subseteq \ilm$, consider the family of complements of sets in $\L$, viz.,  ${\L}^{\ast}:=\{[m]\setminus A : A\in \L\} \subseteq I_{m-\ell}[m]$. Then 
$\L$ is a subclose family in $\ilm$ if and only if $\L^\ast$ is a subclose family in $I_{m-\ell}[m]$. Moreover, 
$$
K_r(\ell,m) =  {{r}\choose {2}}(2\ell-m) + K_r(m-\ell, m) \quad \mbox{ for } 0\le r\le k. 
$$
\end{proposition}

\noindent
\emph{Proof:}
Write $A^c$ for $[m]\setminus A $ for $A\in \ilm$. Then for any $A,B\in \ilm$, we have 
$$
|A \cap B| = m-|A^c\cup B^c| = m-(m-\ell) - (m-\ell) + |A^c\cap B^c| = (2\ell -m) + |A^c\cap B^c|.
$$
Thus if $r:=|\L|$ and we write $\L = \{A_1, \dots , A_r\}$, then 
$$
K_\L = \sum_{1\le i<j\le r}   (2\ell -m) + |A_i^c\cap A_j^c| = {{r}\choose {2}}(2\ell-m) + K_{\L^\ast}.
$$ 
Now as $\L$ varies over families in $\ilm$ of cardinality $r$, the dual $\L^\ast$ varies over families in $I_{m-\ell}[m]$ of cardinality $r$. It follows that
$\L$ is a subclose family in $\ilm$ if and only if $\L^\ast$ is a subclose family in $I_{m-\ell}[m]$. Moreover, 
\begin{equation*}
K_r(\ell,m) =  {{r}\choose {2}}(2\ell-m) + K_r(m-\ell, m) \quad \mbox{ for } 0\le r\le k . 
\eqno{\Box}
\end{equation*}

Recall  that given any $a,b\in \Z$, the 
binomial coefficient ${{a}\choose {b}}$ is defined by 
$$
{{a}\choose {b}}: = \left\{ \begin{array}{ll} \displaystyle{\frac{a(a-1)\cdots (a-b+1)}{b!}} & \ \mbox{ if } b\ge 0, \\ 0 & \ \mbox{ if } b<0.\end{array} \right. 
$$
With this in view, we may permit $a$ and $b$ to take negative values. We 
record 
some elementary properties of binomial coefficients, which 
will be useful in the sequel.  

\begin{lemma}
\label{binomial}
Given any integers $a,b,c,d,e$, we have the following. 
\begin{enumerate}
	\item[{\rm (i)}] $\dbinom{a}{b}=\dbinom{a}{a-b}$ if and only if either $a\geq 0$ or $a<b<0$ 
	\\
	
	\item[{\rm (ii)}] $\dbinom{a}{b}=0$ if and only if either $b<0$ or $b>a\geq 0$
	\\
	
	\item[{\rm (iii)}] $\dbinom{a}{b}\dbinom{b}{c}=\dbinom{a}{c}\dbinom{a-c}{b-c}$. 
	\\
	
	\item[{\rm (iv)}] $\dbinom{a+ b}{c-e}=\displaystyle\sum_{j=e}^c 
	\dbinom{a + d}{c-j}\dbinom{b - d}{j-e}$.
\end{enumerate}
\end{lemma}

\bprf
Both (i) and (ii) are straightforward. Proofs of (iii) and (iv) are also elementary. See, for example, Lemma 3.2 and Corollary 3.4 of \cite{Gh}.
\eprf

The value of $K_r(\ell,m)$ for the maximum permissible parameter $r$ is determined below. 

\begin{proposition}
\label{pro:Kk}
$$
K_k(\ell,m) =m\binom{\nu}{2}, \quad \mbox{ where } \quad \nu:=\binom{m-1}{\ell-1}.
$$ 
\end{proposition}

\bprf
Observe that 
$\frac{m}{\ell}{\nu}=\binom{m}{\ell}=k$, 
i.e., ${\nu}=\ell k/m$.
Write $\ilm =\{A_1, \dots , A_k\}$. Then 
\begin{equation}
\label{U}
 K_{k}(\ell,m) =\displaystyle\sum_{1\leq i<j \leq k}|{A_i}\cap {A_j}|=\frac{1}{2}\left(\sum_{i,j}|{A_i}\cap{A_j}|-
\sum_{1\leq i=j\leq k}\ell\right)=\frac{1}{2}\left[U-k\ell\right], 
\end{equation}
where 
\begin{equation*}
 U:=\displaystyle{ \sum_{A,B\in \ilm}|A\cap B| }.
\end{equation*}
Thus it suffices to determine $U$, which is more symmetric than $K_k(\ell, m)$. To 
find $U$, note that for any $A, B \in \ilm$, the intersection $A\cap B$ is a subset $E$, say,  of $[m]$ of
cardinality $i\leq \ell$. Thus,
\begin{eqnarray*}
U&=&\displaystyle{\mathop{\sum_{E\subseteq[m]}}_{|E|\leq \ell}}\;
\displaystyle{\mathop{\sum_{A,B\in \ilm}}_{A\cap B=E}}|E|\\
&=&\sum_{i=0}^{\ell} {\mathop{\sum_{E\subseteq[m]}}_{|E|=i}} i |\{(A,B) \in \ilm \times \ilm :
A \cap  B = E\}|\\
&=&\sum_{i=0}^{\ell}\dbinom{m}{i}\left[i\dbinom{m-i}{\ell-i}\dbinom{m-\ell
}{\ell-i}\right]\\
&=&\sum_{i=0}^{\ell}i\dbinom{m}{i}\left[\dbinom{m-i}{m-\ell}\dbinom{m-\ell
}{\ell-i}\right]  \qquad \mbox{[by Lemma \ref{binomial} (i)]}\\
&=&\sum_{i=1}^{\ell} m \left[\dbinom{m-1}{m-i}\dbinom{m-i}{m-\ell}\right]\dbinom{m-\ell}{\ell-i} \qquad \mbox{[by Lemma \ref{binomial} (ii)]} \\
&=&m\dbinom{m-1}{m-\ell}\sum_{i=1}^{\ell}\dbinom{\ell-1}{\ell-i}\dbinom{m-\ell}
{\ell-i} \qquad \mbox{[by Lemma \ref{binomial} (iii)]} \\
&=&m\dbinom{m-1}{\ell-1}\sum_{i=1}^{\ell}\dbinom{\ell-1}{i-1}\dbinom{m-\ell}
{\ell-i} \qquad \mbox{[by Lemma \ref{binomial} (i)]}\\
&=&m\dbinom{m-1}{\ell-1}\dbinom{m-1}{\ell-1} \qquad \mbox{[by Lemma \ref{binomial} (iv)]} \\
&= &m{\nu}^2. \\
\end{eqnarray*}
Therefore, equation (\ref{U}) becomes
\begin{eqnarray*} 
K_k(\ell, m)=\frac{1}{2}\left[U-k\ell \right]=\frac{1}{2}\left[m{\nu}^2-m{\nu} \right]=
m\dbinom{{\nu}}{2},
\end{eqnarray*}
as desired. 
\eprf

The above result will help us establish yet another version of duality among subclose families. But first we need a preliminary result whose proof is similar in spirit to the proof above. 

\begin{lemma}
\label{AneB}
Given any $A\in \ilm $, we have
$$
\displaystyle{\mathop{\sum_{ B \in \ilm }}_{B \ne A}} |A\cap B|=\ell({\nu}-1), \quad \mbox{ where }\quad \nu := \binom{m-1}{\ell-1}.
$$
\end{lemma}

\begin{proof}
As in the proof of Proposition \ref{pro:Kk}, we have 
\begin{eqnarray*}
{\mathop{\sum_{ B \in \ilm }}_{B \ne A}} |A\cap B|
&=& \sum_{i=0}^{\ell-1}\; {\mathop{\sum_{E\subseteq A}}_{|E|=i}}\;{\mathop{\sum_{ B \in \ilm}}_{ B \cap A =E}} |E| 
\\
&=& \sum_{i=0}^{\ell-1} \binom{\ell}{i} i \binom{m-\ell}{\ell-i} \\
&=&\ell\sum_{i=1}^{\ell-1}\binom{\ell-1}{i-1} \binom{m-\ell}{\ell-i}\\
&=&\ell\left[-1 + \sum_{i=1}^{\ell}\binom{\ell-1}{i-1}\binom{m-\ell}{\ell-i}\right]\\
&=&\ell\left[\binom{m-1}{\ell-1}-1\right], 
\end{eqnarray*}
where the last equality follows from part (iv) of Lemma \ref{binomial}.
\end{proof}

\begin{proposition}[Second Duality Theorem]
\label{secondduality} 
Given $\L\subseteq \ilm$, consider the complement  ${\L}^{c}:=\ilm \setminus \L$. 
Then $\L$ is a subclose family in $\ilm$ if and only if $\L^c$ is a subclose family in $I_{\ell}[m]$. Moreover, if we let $r:=|\L|$ and $\nu := \binom{m-1}{\ell-1}$, 
then 
\begin{equation}
\label{KLc}
K_{\L^c} = m {{\nu}\choose {2}} - r\ell(\nu -1) + K_{\L}.
\end{equation}
Consequently, 
\begin{equation}
\label{Kkminusr}
K_{k-r}(\ell,m) =  m {{\nu}\choose {2}} - r\ell(\nu -1) + K_{r}(\ell, m) \quad \mbox{ for } 0\le r\le k. 
\end{equation}
\end{proposition}

\bprf
Let $\L \subseteq \ilm$. In view of Proposition \ref{pro:Kk}, we have 
$$
m {{\nu}\choose {2}} = K_k(\ell,m) = K_{\ilm} = K_\L + K_{\L^c} + \sum_{A\in \L}\sum_{B\in \L^c} |A\cap B|.
$$
Further, in view of Lemma \ref{AneB}, we can write 
\begin{eqnarray*}
\sum_{A\in \L}\sum_{B\in \L^c} |A\cap B| 
&=& \sum_{A\in \L} \left( {\mathop{\sum_{ B \in \ilm}}_{ B \ne A}} |A\cap B| 
- {\mathop{\sum_{ B \in \L}}_{ B \ne A } |A\cap B|} \right) \\
&=& \sum_{A\in \L} \ell(\nu -1) - \sum_{A\in \L}{\mathop{\sum_{ B \in \L}}_{ B \ne A } |A\cap B|}\\
&=& r \ell(\nu -1) - 2 K_{\L}.
\end{eqnarray*}
It follows that 
$$
K_{\L^c} = m {{\nu}\choose {2}} - r\ell(\nu -1) + K_{\L}.
$$
This implies that $K_{\L^c} \le m {{\nu}\choose {2}} - r\ell(\nu -1) + K_{r}(\ell,m)$, and the equality holds if and only if $\L$ is subclose. Consequently, 
$$
K_{k-r}(\ell,m) = m {{\nu}\choose {2}} - r\ell(\nu -1) + K_r(\ell,m),
$$
and moreover, $\L$ is subclose if and only if $\L^c$ is subclose. 
\eprf

\begin{corollary}
\label{dualuptomu}
If $s\in \Z$ is such that $k-\mu \le s \le k$, then 
$$
K_{s}(\ell,m) = m {{\nu}\choose {2}} - \ell(\nu -1)(k-s) + (\ell -1) {{k-s}\choose {2}}, \quad \mbox{ where }\quad \nu := \binom{m-1}{\ell-1}.
$$
In particular, if $m\ge 4$, then 
$$
K_s(2,m) = m {{m-1}\choose {2}} - 2\left(k-s\right)(m -2) +  {{k-s}\choose {2}}, \quad \mbox{ for } {{m-1}\choose{2}} \le s \le {{m}\choose{2}}.
$$
\end{corollary}

\bprf
Given $s\in \Z$ with $k-\mu \le s \le k$, observe that $r:=k-s$ satisfies $0\le r\le \mu$. Now use \eqref{Kkminusr} 
together with Proposition \ref{subcloseisclose}
to obtain the first equality. The second equality follows from the first by noting that if $\ell =2$ and $m\ge 4$, then $\mu=m-1=\nu$   and $k-\mu = {{m-1}\choose{2}}$.  
\eprf

\begin{example}
Using the above results, one can readily compile a table of values of $K_r(\ell,m)$ for $\ell=2$ and for small values of $m$. For example, we have

\medskip

\begin{center}
\begin{tabular}{|c|c|c|c|c|c|c|c|c|c|c|}
\hline 
$r$ & 1 & 2 & 3 & 4 &5 & 6 & 7 & 8 & 9 & 10 \\ \hline 
$K_r(2,5)$ & 0 & 1 & 3 & 6 & 8 & 12 & 15 & 19 & 24 & 30 \\ \hline
\end{tabular}
\end{center}
\bigskip
\noindent
and

\bigskip
\begin{center}
\begin{tabular}{|c|c|c|c|c|c|c|c|c|c|c|c|c|c|c|c|}
\hline 
$r$ & 1 & 2 & 3 & 4 &5 & 6 & 7 & 8 & 9 & 10 & 11 & 12 & 13 & 14 & 15\\ \hline 
$K_r(2,6)$ & 0 & 1 & 3 & 6 & 10 & 12 & 15 & 19 & 24 & 30  & 34 & 39 & 45 & 52 & 60\\ \hline
\end{tabular}

\bigskip

\end{center}

\noindent
where it may be noted that the barrier $\mu$ on the values of $r$ is $4$ in the first table and $5$ in the second table. That is where the pattern begins to change. 
\end{example}

\section{Higher Weights of Grassmann Codes and Schubert Codes}
\label{sec2}

We have made it amply clear in the Introduction that the combinatorial considerations in the preceding section were motivated by problems in Coding Theory, more specifically, the determination of the higher weights of linear codes associated to Grassmann and Schubert varieties. In this section, we begin by describing some relevant background, set up some notation, and then state a precise conjecture relating the said higher weights to subclose families. 

Fix integers $k,n$ with $1\le k \le n$ and a prime power $q$. 
Let $C$ be a linear $[n,k]_q$-code, i.e., 
let $C$ be a $k$-dimensional subspace of the $n$-dimensional vector space $\Fq^n$ 
over the finite field $\Fq$ with $q$ elements. 
Given any $x=(x_1,\dots,x_n)$ in $\Fq^n$, let
$$
\supp (x):=\{ i : x_i\ne 0\} \quad \mbox{ and } \quad \Vert x \Vert := |\supp (x)|
$$
denote the \emph{support} and the \emph{(Hamming) norm} of $x$. 
More generally, for $D\subseteq \Fq^n$, 
let
$$
\supp (D):=\{ i : x_i\ne 0 \mbox{ for some } x=(x_1,\dots ,x_n) \in D\} 
\quad \mbox{ and } \quad \Vert D \Vert := |\supp (D)|
$$
denote the \emph{support} and the \emph{(Hamming) norm} of $D$. The \emph{minimum distance} or the \emph{Hamming weight} of $C$ is defined by 
$d(C):=\min\{\Vert x \Vert : x\in C \mbox{  with } x\ne 0\}$. 
More generally, following \cite{W}, for any positive integer $r$, the $r^{\text{th}}$ 
\emph{higher weight}  or the $r^{\text{th}}$ 
\emph{generalized Hamming weight} $d_r = d_r(C)$ of the code $C$ is defined by
$$
d_r(C) := \min\left\{ \Vert D \Vert : D \mbox{ is a subspace of $C$ with } \dim D =r\right\}.
$$
Note that $d_1(C)=d(C)$. 
If $C$ is \emph{nondegenerate}, i.e., if
$C$ is not contained in a coordinate hyperplane of $\Fq^n$, then it is easy to see that
$$ 
0< d_1(C) < d_2(C) < \cdots < d_k(C) = n.
$$
See, for example, \cite{TV2} for a proof as well as a great deal of basic information about
higher weights of codes. The set $\{d_r(C):  1\le r \le k\}$ is often referred to as the
\emph{(complete) weight hierarchy} of the code $C$. It is usually interesting, and difficult, to   
determine the complete weight hierarchy of a given code. 

An equivalent way of describing codes is via the language of projective systems, explained, for example in \cite{TV2,N,GL}. A $[n,k]_q$-projective system $X$ is a (multi)set of $n$ points in the projective space $\PP^{k-1}$ over $\Fq$. We say $X$ is \emph{nondegenerate} if it is not contained in a hyperplane of $\PP^{k-1}$. An $[n,k]_q$-nondegenerate projective system gives rise to a unique (up to isomorphism) nondegenerate $[n,k]_q$-linear code $C_X$. The minimum distance of $C_X$ corresponds to maximizing the number of points of hyperplane sections of $X$, while the $r^{\text{th}}$ higher weight corresponds to  maximizing the number of points of  sections of $X$ by codimension $r$ projective linear subspaces. More precisely, for $0\le r \le k$, we have 
$$
d_r(C_X) = 
n - \max \left\{ |X \cap {\Pi}| : {\Pi} \mbox{ is a projective 
subspace of codimension $r$ in } \PP^{k-1} \right\} .
$$

Linear codes associated to projective systems given by the $\Fq$-rational points of higher dimensional projective algebraic varieties defined over $\Fq$ have been of much interest lately, and we refer to the recent survey by  
Little \cite{little} for more on this. 
We are particularly interested in the case of 
Grassmann variety $G_{\ell,m}$ and its Schubert subvarieties $\Omega_{\alpha}= \Omega_{\alpha}(\ell,m)$ with its nondegenerate Pl\"ucker embedding in $\PP^{k-1}$ and $\PP^{k_{\alpha}-1}$, respectively. Here, as in Section \ref{sec1}, $\ell,m$ are fixed positive integers with $\ell \le m$ and $k:={{m}\choose{\ell}}$, while   
$\alpha$ varies over the natural indexing set for points of $\PP^{k-1}$, namely, 
$$
I(\ell,m):=\{ \b = (\b_1, \dots , \b_{\ell} )\in \Z^{\ell}  : 
1\le \b_1 < \dots < \b_\ell  \le m \},
$$
and for any $\a\in I(\ell,m)$, 
$$
k_{\a}: = \left|I_{\a}(\ell,m)\right| \quad \mbox{where} \quad 
I_{\a}(\ell,m):=\{ \b \in I(\ell,m) : \b_i \le \a_i \mbox{ for all } i=1,  \dots , \ell \}.
$$
We identify 
$\PP^{k_{\alpha}-1}$ 
with $\left\{p=(p_{\b})\in \PP^{k-1} : p_{\b} = 0 \mbox{ for all } \beta \in I(\ell,m)\setminus I_{\a}(\ell,m) \right\}$ so that $\Omega_{\alpha}(\ell,m) = G_{\ell,m}\cap \PP^{k_{\alpha}-1}$. 
For precise definitions of $G_{\ell,m}$ and $\Omega_{\alpha}$, 
and their Pl\"ucker embeddings, we refer to \cite{GL} and \cite{GT} or the references therein. The linear codes corresponding to $G_{\ell,m}$ and $\Omega_{\alpha}(\ell,m)$ are denoted by $C(\ell,m)$ and $C_{\alpha}(\ell,m)$ respectively. 
The length $n$ of  $C(\ell,m)$ and $n_{\alpha}$ of $C_{\alpha}(\ell,m)$ are respectively given by 
$$
n = |G_{\ell,m}(\Fq)| = { {m} \brack {\ell} }_q :=
 \frac{ (q^m -1)(q^{m} - q) \cdots (q^{m} - q^{\ell -1})}
{(q^{\ell} -1)(q^{\ell} - q) \cdots (q^{\ell} - q^{\ell -1})} \ \mbox{ and } \  
n_{\a}= |\Omega_{\alpha}(\Fq)|.
$$
The dimension $k$ of $C(\ell,m)$ and $k_{\alpha}$ of $C_{\alpha}(\ell,m)$ are respectively given by 
$$
k = |I(\ell,m)| = {{m}\choose{\ell}} \quad \mbox{and} \quad k_{\a}: = \left|I_{\a}(\ell,m)\right|.
$$
A number of explicit formulas for 
$n_{\a}$ and $k_{\a}$ are given in \cite{GT}.

Given any  $\Lambda \subseteq I(\ell,m)$, we let $\Pi_{\L}$ denote the intersection of the 
corresponding Pl\"ucker coordinate hyperplanes; more precisely, 
$$
\Pi_{\L} : = \left\{p=(p_{\b})\in \PP^{k-1} : p_{\b} = 0 \mbox{ for all } \beta \in \Lambda \right\}.
$$
Note that $\Pi_{\L}$ is a  projective linear subspace of codimension $|\L|$ in $\PP^{k-1}$, and also that if 
$\Lambda \subseteq I_{\a}(\ell,m)$, then $\Pi_{\L}\cap \PP^{k_{\alpha}-1}$ is a projective linear subspace of 
codimension $|\L|$ in $\PP^{k_{\alpha}-1}$. 

There is a natural one-to-one correspondence between the indexing set $I(\ell,m)$ and the set $I_{\ell}[m]$ defined in the previous section, 
given simply by 
$$
\b = (\b_1, \dots , \b_{\ell} ) \longmapsto \bar{\b}=\{\b_1, \dots , \b_{\ell} \}.
$$
With this in view, we shall identify $I(\ell,m)$ with $I_{\ell}[m]$, and apply the notions and results of Section \ref{sec1} for $I_{\ell}[m]$ and its subfamilies to $I(\ell,m)$ and its subfamilies. In particular, we can talk about subclose families in $I(\ell,m)$. 
We are now ready to propose a plausible fact  about the higher weights of  $C(\ell,m)$ and $C_{\alpha}(\ell,m)$.

\begin{conjecture}
Let $r$ be a positive integer. If $\, r\le k$, then 
the $r^{\text{th}}$ higher weight  of the Grassmann code  $C(\ell,m)$ is given by 
$$
d_r(C(\ell,m)) = 
{{m} \brack {\ell} }_q - \; \max \left\{ |G_{\ell,m}(\Fq) \cap \Pi_{\Lambda}| : {\L}\subseteq  I(\ell,m) \mbox{ is  subclose and } |\L|=r\right\}.
$$
More generally, given any $\a\in I(\ell,m)$, if $\, r\le k_{\a}$,  then the $r^{\text{th}}$ higher weight of the Schubert code $C_{\a}(\ell,m)$ is given by
$$
d_r(C_{\a}(\ell,m)) = 
n_{\a} - \max \left\{ |\Omega_{\alpha}(\Fq) \cap \Pi_{\Lambda}| : {\L}\subseteq  I_{\a}(\ell,m) \mbox{ is  subclose and } |\L|=r\right\}.
$$
\end{conjecture}

The evidence we have in favor of this conjecture is as follows. 

\begin{enumerate}
	\item 
	The conjecture is true in the case of $C(\ell,m)$ for $1\le r \le \max\{\ell, m-\ell\} + 1$. (See \cite{GL}.)
	\item
	The conjecture is true in the case of $C(2,m)$ for $r= \max\{2, m-2\} + 2$. (See \cite{GPP}.)
	\item
	The conjecture is true in the case of $C_{\a}(\ell,m)$ for $r=1$. (See \cite{xiang}.)
	\item
	The conjecture is true in the case of $C_{\a}(\ell,m)$ where $\a$ is a submaximal element of $I_{\a}(\ell,m)$ [so that the corresponding Schubert variety $\Omega_{\a}$ is of codimension $1$ in $G_{\ell,m}$] for 
	$1\le r \le \max\{\ell, m-\ell\} $. (See \cite{GT}.)
\end{enumerate}
It may be remarked that the notion of a subclose family and the above conjecture is similar to, yet distinct from, the notion of a Schubert union introduced in \cite{HJR} and the corresponding conjecture of Hansen, Johnsen and Ranestad \cite{HJR,HJR2} that the higher weights are attained by Schubert unions. It may also be noted that for a given $r$, there may be more than one subclose family of cardinality $r$. Thus the above conjecture does not pinpoint to a single such family but simply narrows down the search for such a family. 

\section{Threshold Graphs, Optimal Graphs and Subclose Families}
\label{sec3}

When $\ell=2$, the elements of the family $\ilm$ of $\ell$-subsets of $[m]:=\{1,\dots ,m\}$ can be viewed as the edges of a graph. It is, therefore, natural to investigate if the combinatorial notions and results in Section \ref{sec1} have analogues and extensions in the rich and diverse field of graph theory. We will attempt to address these concerns in this section. 

Given any $\L\subseteq I_2[m]$, we denote by $G_{\L}$ the graph whose vertex set is $[m]$ and the edge set is $\L$. Note that this is a simple (undirected) graph. Conversely, any simple graph on $[m]$ is of the form $G_{\L}$ for a unique  $\L\subseteq I_2[m]$. To say that $\L$ is close corresponds to saying that any two edges of $G_{\L}$ are incident. 
Thus, Proposition \ref{StrThm} corresponds to the following elementary result in graph theory.

\begin{proposition} 
\label{StrThmGraphs}
A simple graph in which any two edges are incident is either a star or a triangle. 
\end{proposition}

The analogue of subclose family is more interesting. Before explaining this, let us recall some notions from graph theory.\footnote{We are using here a notation that is consistent with the notation of Section \ref{sec1}. Inconvenience caused, if any,  to graph theorists, who may be more used to letting $n$ be the number of vertices, $e$ the number of edges, and $d_i$ the degree of the vertex $i$, is regretted.}
Let $G$ be a $(m,r)$-graph, i.e., a graph with $m$ vertices (assumed to be elements of the set $[m]$) and $r$ edges. We denote by $g_i=g_i(G)$ the \emph{degree} of the vertex $i$, viz., the number of edges emanating from it. The sequence $(g_1, \dots , g_m)$ is called the \emph{degree sequence} of $G$. It is well-known and easy to see that
\begin{equation}
\label{sumdeg}
\sum_{i=1}^m g_i = 2r.
\end{equation}
A simple graph $G$ is said to be a \emph{threshold graph} if $G$ can be constructed from a one-vertex graph by repeatedly adding an isolated vertex or a universal one (i.e., a vertex adjacent to every other vertex). 
A simple $(m,r)$-graph  $G$ is said to be $(m,r)$-\emph{optimal}, or simply, \emph{optimal} if 
$$
\Sigma(G):= \sum_{i=1}^m g_i(G)^2
$$
is maximum among all simple $(m,r)$-graphs. Threshold graphs are a topic of considerable interest in graph theory, and we refer to \cite{MP} for more on this. It is easy to see that an optimal graph is a threshold graph (cf. \cite[Fact 3]{peled}). In \cite{peled} it is shown that an optimal graph is one among certain six explicit classes of graphs. 
However, as the authors of \cite{peled} say, the complete characterization of optimal graphs remains an open question. 
The following explicit bound for $\Sigma(G)$ for a $(m,r)$-graph $G$ is given by de Caen \cite{caen}.
\begin{equation}
\label{sigmag}
\Sigma(G) \le C(r,m) \quad \mbox{ for } m\ge 2, \quad \mbox{ where } \quad 
C(r,m) : = r \left(\frac{2r}{m-1} + m-2\right).
\end{equation}
A somewhat more general bound has been obtained by Das \cite{das}; however, this bound is not a pure function of $m$ and $r$, but involves the maximum and the minimum among the vertex degrees $g_1, \dots , g_m$. 

The relation between optimal graphs and subclose families is given below. 

\begin{proposition}
\label{subcloseoptimal}
Let $\L$ be a subset of $I_2[m]$ with $|\L|=r$. Then
$$ 
K_{\L} = \frac 12 \Sigma(G_{\L}) - r \quad \mbox{or equivalently,} \quad \Sigma(G_{\L}) = 2K_{\L} + 2r.
$$
Consequently, $\L$ is a subclose family if and only if $G_{\L}$ is an optimal graph. Moreover,
$$
\max\{\Sigma(G): G \mbox{ is a simple $(m,r)$-graph}\} = 2K_r(2,m) + 2r.
$$
\end{proposition}

\bprf
Write $\L=\{A_1, \dots , A_r\}$. Note that $A_i\cap A_j$ is either empty or singleton for $1\le i<j\le r$. Thus,  
$$
K_{\L} = \sum_{i<j} |A_i\cap A_j| = \sum_{i<j} \; \sum_{v\in A_i\cap A_j} 1 = \sum_{v\in [m]} {\mathop{\sum_{i<j}}_{v\in A_i\cap A_j}} 1 =  \sum_{v\in [m]} {{g_v}\choose {2}}.
$$
Further, in view of \eqref{sumdeg}, we have 
$$
K_{\L} = \frac 12 \sum_{v\in [m]} g_v^2 - r = \frac 12 \Sigma(G_{\L}) - r.
$$
Since $G_{\L}$ varies over all  simple $(m,r)$-graphs as $\L$ varies over subsets of $I_2[m]$ of cardinality $r$, it follows that $\L$ is subclose if and only if $G_{\L}$ optimal, and moreover,
$$
\max\{\Sigma(G): G \mbox{ is a simple $(m,r)$-graph}\} = 2K_r(2,m) + 2r,
$$
as desired. 
\eprf

\begin{corollary} Assume that $m\ge 4$. Then for any  simple $(m,r)$-graph $G$, we have 
\begin{equation}
\label{trivialbound}
\Sigma(G) \le r(r+1) \quad \mbox{ for } r\le m-1, 
\end{equation}
and the equality holds if and only if $G$ is a star with $r+1$ vertices (and $m-r-1$ isolated vertices). 
\end{corollary}

\bprf
Since $m\ge 4$, we have $\mu:=\max\{2,m-2\}+1 = m-1$. Thus, $K_r(2,m)\le {{r}\choose {2}}$ for $r\le m-1$, thanks to Proposition \ref{subcloseisclose}. Now apply Proposition \ref{subcloseoptimal} to obtain \eqref{trivialbound}.
The assertion about the equality follows from Proposition \ref{StrThmGraphs}. 
\eprf
Already, the trivial bound given by the Corollary above is superior to de Caen's bound
\eqref{sigmag} in several cases. Indeed, an easy calculation shows that
$$
C(r,m) - r(r+1) = \frac{r(m-3)(m-1-r)}{m-1} > 0 \quad \mbox{ for } r < m-1 \mbox{ and } m\ge 4.
$$
The above result may also be compared with one of the cases where equality holds in a bound for $\Sigma (G)$ obtained by Das \cite{das}.

It is well-known that dual graphs $\overline{G}$ of optimal graphs $G$ [defined in such a way that the non-edges of $G$ are the edges of $\overline{G}$] are optimal (cf. \cite[Fact 1]{peled}). This corresponds precisely to the special case $\ell=2$ of the Second Duality Theorem (Proposition \ref{secondduality}). Further, it is easy to  see that when $\ell=2$, equation \eqref{KLc} corresponds precisely to the following elementary relation for simple $(m,r)$-graphs $G$:
$$
\Sigma(\overline{G}) = m(m-1)^2 - 4r(m-1) + \Sigma(G).
$$
The following bound, dual to the trivial bound given by \eqref{trivialbound}, appears to be new. 

\begin{corollary}
If $G$ is a simple $(m,r)$-graph such that ${{m-1}\choose{2}} \le r \le {{m}\choose{2}}$, then 
$$
\Sigma(G)\le m(m-1)(m-2) + (k-r)(k-r-1) - 4(k-r)(m-2)+2r, \quad \mbox{where} \quad k:={{m}\choose{2}}.
$$
\end{corollary}

\bprf
Follows from Corollary \ref{dualuptomu} and Proposition \ref{subcloseoptimal}.
\eprf

Finally, we remark that the First Duality Theorem (Proposition \ref{firstduality}) has no analogue in the setting of optimal graphs for the simple reason that it relates optimal graphs to objects that are not graphs, but hypergraphs. Indeed, as far as we know, not much seems to be known about threshold hypergraphs and optimal hypergraphs. Perhaps the notion of a subclose family and the results of Section \ref{sec1} may be of some help in this direction. 

\section*{Acknowledgments}

We are grateful to Murali Srinivasan for helpful discussions and bringing \cite{peled} to our attention.

\end{document}